# EXIT PROBLEM OF A TWO-DIMENSIONAL RISK PROCESS FROM THE QUADRANT: EXACT AND ASYMPTOTIC RESULTS

By Florin Avram,[1,2] Zbigniew Palmowski[2,3]
and Martijn R. Pistorius[1,4]

*Université de Pau, University of Wrocław and King's College London*

Consider two insurance companies (or two branches of the same company) that divide between them both claims and premia in some specified proportions. We model the occurrence of claims according to a renewal process. One ruin problem considered is that of the corresponding two-dimensional risk process first leaving the positive quadrant; another is that of entering the negative quadrant. When the claims arrive according to a Poisson process, we obtain a closed form expression for the ultimate ruin probability. In the general case, we analyze the asymptotics of the ruin probability when the initial reserves of both companies tend to infinity under a Cramér light-tail assumption on the claim size distribution.

**1. Introduction.**

*The multidimensional renewal risk model.* In collective risk theory, the reserves process $X$ of an insurance company is modeled as

$$X(t) = x + pt - S(t), \qquad (1)$$

where $x$ denotes the initial reserve, $p$ is the premium rate per unit of time and $S(t)$ is a stochastic process modeling the amount of cumulative claims up to time $t$. Taking $S$ to be a compound Poisson or com-

Received June 2006; revised February 2008.
[1]Supported by London Mathematical Society Grant 4416.
[2]Supported by POLONIUM no 09158SD.
[3]Supported by the Ministry of Science and Higher Education of Poland Grant N N2014079 33 (2007–2009) and NWO 613.000.310.
[4]Supported by Nuffield Foundation Grant NUF/NAL/000761/G and EPSRC Grant EP/D039053/1.

pound renewal process yields the Cramér–Lundberg model and the Sparre–Andersen model, respectively. Recently, several authors have studied extensions of classical risk theory toward a multidimensional reserves model (1) where $X(t), x, p$ and $S(t)$ are vectors, with possible dependence between the components of $S(t)$. Indeed, the assumption of independence of risks may easily fail, for example, in the case of reinsurance, when incoming claims have an impact on both insuring companies at the same time. In general, one can also consider situations where each claim event might induce more than one type of claim in an umbrella policy [see Sundt (1999)]. For some recent papers considering dependent risks, see Dhaene and Goovaerts (1996, 1997), Goovaerts and Dhaene (1996), Müller (1997a, 1997b) and Denuit, Genest and Marceau (1999), Ambagaspitiya (1999), Dhaene and Denuit (1999), Hu and Wu (1999) and Chan, Yang and Zhang (2003).

*Model and problem.* In this paper we consider a particular two-dimensional risk model in which two companies split the amount they pay out of each claim in proportions $\delta_1$ and $\delta_2$ where $\delta_1 + \delta_2 = 1$, and receive premiums at rates $c_1$ and $c_2$, respectively. Let $U_i$ denote the risk process of the $i$th company

$$U_i(t) := -\delta_i S(t) + c_i t + u_i, \qquad i = 1, 2,$$

where $u_i$ denotes the initial reserve. We will study here the eventual ruin probabilities in two cases:

1. the Lévy model, obtained by taking $S(t)$ to be a general Lévy process;
2. the Sparre–Andersen/renewal risk model, where $S(t)$ is

$$(2) \qquad S(t) = \sum_{i=1}^{N(t)} \sigma_i,$$

$N(t)$ is a renewal process with i.i.d. interarrival times $\zeta_i$, and the claims $\sigma_i$ are i.i.d. nonnegative random variables independent of $N(t)$.

The intersection of the two cases is the classical Cramér–Lundberg model, where $S(t)$ is a compound Poisson process with nonnegative jumps.

We shall denote by $F(x)$ the distribution function of the "claims" $\sigma_i$, and by $\lambda$ and $\mu$ the reciprocals of the means of $\zeta_i$ and $\sigma_i$, respectively.

As usual in risk theory, we assume that $U_i(t) \to \infty$ a.s. as $t \to \infty$ ($i = 1, 2$). In the case of the Sparre–Andersen model, this amounts to $p_i > \rho := \frac{\lambda}{\mu} = E\sigma/E\zeta$.

We shall assume that the second company, to be called reinsurer, receives less premium per amount paid out, that is,

$$(3) \qquad p_1 = \frac{c_1}{\delta_1} > \frac{c_2}{\delta_2} = p_2.$$

Several ruin problems will be considered here:

1. The first time $\tau_{\mathrm{or}}$ when (at least) one insurance company is ruined, that is, the exit time of $(U_1(t), U_2(t))$ from the positive quadrant

(4) $$\tau_{\mathrm{or}}(u_1, u_2) := \inf\{t \geq 0 : U_1(t) < 0 \text{ or } U_2(t) < 0\}.$$

2. The first time $\tau_{\mathrm{sim}}$ when the insurance companies experience simultaneous ruin, that is, the entrance time of $(U_1(t), U_2(t))$ into the negative quadrant

(5) $$\tau_{\mathrm{sim}}(u_1, u_2) := \inf\{t \geq 0 : U_1(t) < 0 \text{ and } U_2(t) < 0\}.$$

The associated ultimate/perpetual ruin probabilities will be respectively denoted by $\psi_{\mathrm{or}}(u_1, u_2)$ and $\psi_{\mathrm{sim}}(u_1, u_2)$

$$\psi_{\mathrm{or}}(u_1, u_2) = P(\tau_{\mathrm{or}}(u_1, u_2) < \infty),$$
$$\psi_{\mathrm{sim}}(u_1, u_2) = P(\tau_{\mathrm{sim}}(u_1, u_2) < \infty).$$

Letting $\tau_i(u_i) = \inf\{t \geq 0 : U_i(t) < 0\}$, $i = 1, 2$, we also will consider

(6) $$\psi_{\mathrm{and}}(u_1, u_2) = P(\tau_1(u_1) < \infty \text{ and } \tau_2(u_2) < \infty).$$

Denoting by $\psi_i(u) := P(\tau_i(u) < \infty)$, the ruin probability of $U_i$ when $U_i(0) = u$, it clearly holds that

(7) $$\psi_{\mathrm{sim}}(u_1, u_2) \leq \psi_{\mathrm{and}}(u_1, u_2) = \psi_1(u_1) + \psi_2(u_2) - \psi_{\mathrm{or}}(u_1, u_2).$$

As the ruin probabilities $\psi_{\mathrm{or}}, \psi_{\mathrm{and}}$ and $\psi_{\mathrm{sim}}$ do not change under a scaling of $(U_1, U_2)$, we restrict ourselves in the sequel to the respective ruin probabilities $\psi_{\mathrm{or}}, \psi_{\mathrm{and}}$ and $\psi_{\mathrm{sim}}$ of the scaled process $(X_1, X_2)$ given by $X_i(t) := U_i(t)/\delta_i = x_i + p_i t - S(t)$ with $x_i = u_i/\delta_i$, $p_i = c_i/\delta_i$. This puts into evidence the fact that all the randomness in our model acts in one direction; in the future, we call this model informally a "two-dimensional degenerate model."

*Geometrical considerations.* The solutions of the "degenerate two-dimensional" ruin problems $\psi_{\mathrm{or}}, \psi_{\mathrm{sim}}$ and $\psi_{\mathrm{and}}$ strongly depend on the relative position of the vector of premium rates $p = (p_1, p_2)$ with respect to the proportions vector $(1, 1)$. A key observation is that the ruin times $\tau_{\mathrm{or}}$ and $\tau_{\mathrm{sim}}$ are also equal to

(8) $$\tau_b(x_1, x_2) = \inf\{t \geq 0 : S(t) > b(t)\},$$

where $b = b_{\min} = \min\{\ell_1, \ell_2\}$ and $b = b_{\max} = \max\{\ell_1, \ell_2\}$, respectively, with

(9) $$\ell_i(t) = x_i + p_i t, \qquad i = 1, 2, t \geq 0;$$

see Figure 1. The two-dimensional problems may thus also be viewed as a one-dimensional crossing problem over a (piecewise) linear barrier.

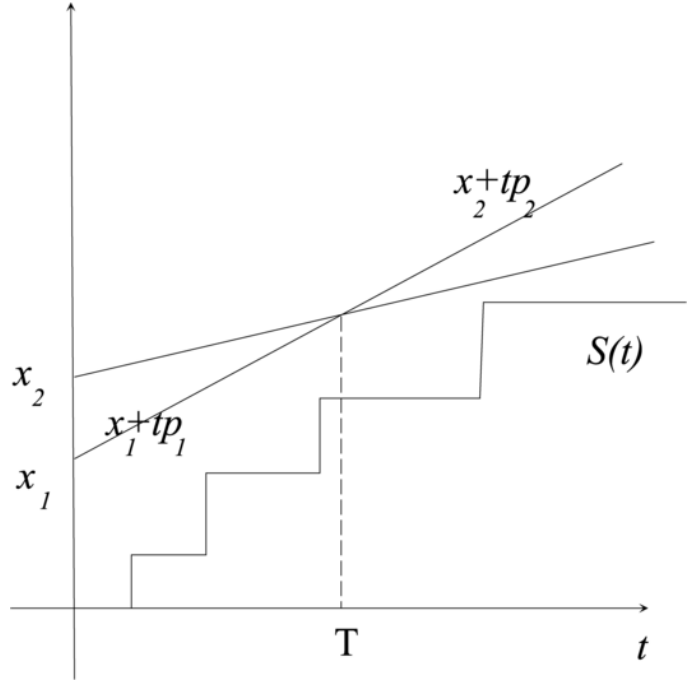

FIG. 1. *The piecewise-linear barrier corresponding to the degenerate two-dimensional first passage problem:* $b_{\min}(t) = \min_{i=1,2}\{x_i + p_i t\}, b_{\max}(t) = \max_{i=1,2}\{x_i + p_i t\}$.

Our exact results follow directly from this geometrical observation, which essentially breaks the problem in two pieces: ruin of one of the coordinates before the deterministic time

$$T = T(x_1, x_2) = \frac{(x_2 - x_1)_+}{p_1 - p_2} \tag{10}$$

of entering the lower cone $x_2 \leq x_1$, or ruin of the other coordinate subsequently.

If the initial reserves satisfy $x_2 \leq x_1$, the two lines do not intersect. It follows therefore that the barriers $b_{\min}, b_{\max}$ are actually linear

$$b_{\min}(t) = x_2 + p_2 t, \qquad b_{\max}(t) = x_1 + p_1 t.$$

Thus, the "or" and "sim" ruin always happen for the second and first company, respectively, and

$$\psi_{\mathrm{or}}(x_1, x_2) = \psi_2(x_2), \qquad \psi_{\mathrm{sim}}(x_1, x_2) = \psi_1(x_1).$$

This case in which explicit and asymptotic formulas for $\psi_{\mathrm{sim}}, \psi_{\mathrm{or}}$ and $\psi_{\mathrm{and}}$ follow directly from the well-known one-dimensional ruin theory; see, for example, Rolski et al. (1999) or Asmussen (2000)—will not be discussed further.

*Contents.* This paper is devoted to obtaining exact results for the eventual ruin probabilities in the other case, when the initial capitals satisfy $x_1 < x_2$. We obtain also sharp "Cramér/light tails" asymptotics, as the initial capitals $x_1, x_2$ tend to infinity along a ray (i.e., $x_1/x_2$ is constant).

We introduce notation and gather some necessary prerequisites from one dimensional ruin theory in Section 2.

In Section 3, we obtain several exact decomposition formulas for the two-dimensional ruin probabilities $\psi_{\text{or}}$, $\psi_{\text{sim}}$ and $\psi_{\text{and}}$, in terms of one-dimensional ruin probabilities. In particular, we obtain for the compound Poisson case with exponential jumps an exact result in Corollary 1.

In Section 4, we show in Theorem 3 for the general two-dimensional Lévy case the natural result:

$$\psi_{\text{or}}(x_1, x_2) \sim \psi_1(x_1) + \psi_2(x_2),$$
$$\psi_{\text{and}}(x_1, x_2) = \text{o}(\psi_1(x_1) + \psi_2(x_2)),$$

where we write $f(x) = \text{o}(h(x))$ $(x \to \infty)$ if $\lim_{x \to \infty} f(x)/h(x) = 0$ and similarly $f(x) \sim h(x)$ $(x \to \infty)$ if $\lim_{x \to \infty} f(x)/h(x) = 1$. The corresponding renewal version is stated in Theorem 4.

In Section 5, specializing to the case of the Cramér–Lundberg model, we sharpen the general result, obtaining two term asymptotic expansions—Theorems 5 and 6. We find different leading terms within subcones of the positive quadrant, as typical in such cases; see, for example, Borovkov and Mogulskii (2001) and Ignatyuk, Malyshev and Shcherbakov (1994). We also find in our specific "degenerate model" a correction term not present in previous works.

The paper concludes with two explicit examples in Section 6, which provides interesting illustrations of multidimensional *sharp* large deviations—see the Appendix for the relation to the existing first passage large deviations theory. In this context, it is worth clarifying that our particular results are considerably sharper than the general logarithmic asymptotics results obtained, for example, by Collamore (1996), Theorems 2.1 and 2.2, for the first passage times to convex open and open sets in $\mathbb{R}^d$, respectively. Further, our proofs do not appeal at all to the multidimensional large deviations theory. Instead, exploiting the special structure of our problem, we make use of a one-dimensional asymptotic limit result; see Theorem 2 in Section 2.3—which is a consequence of the sharp approximation of finite-time ruin probabilities obtained by Arfwedson (1955) and Höglund (1990).

**2. Preliminaries: One dimensional theory.** Let $Z$ be a general Lévy process, that is, a process with stationary and independent increments that is continuous in probability and starts at 0, defined on some probability space $(\Omega, \mathcal{F}, P)$ and let $E$ denote the expectation w.r.t. $P$. To avoid trivialities, we exclude the case that $Z$ has monotone paths.

We will restrict ourselves to Lévy processes $Z$ that admit negative exponential moments, that is, $E[e^{-\nu Z(t)}] < \infty$ for some $\nu > 0$. For such processes, $Z$, we consider the *cumulant exponent*:

$$\kappa(\theta) := t^{-1} \log(E[e^{\theta Z(t)}]),$$

which is well defined on some maximal domain

$$\Theta = \{\theta \in \mathbb{R} : \kappa(\theta) < \infty\}$$

whose interior will be denoted by $\Theta^o := (\underline{\theta}, \overline{\theta})$, where $\underline{\theta} = \inf \Theta$ and $\overline{\theta} = \sup \Theta$. The map $\theta \mapsto \kappa(\theta)$ restricted to the interval $(\underline{\theta}, \overline{\theta})$ is a convex differentiable function. By $\kappa'(\theta)$ $[\kappa'_+(\theta)]$, we denote the [left-]derivative of $\kappa$ at $\theta$, respectively.

In the particular case that $Z$ is equal to the classical Cramér–Lundberg process (i.e., a positive drift added to a spectrally negative compound Poisson process), we note that $\overline{\theta} = \infty$ and $\kappa'(\theta) \to \infty$ if $\theta \to \infty$.

2.1. *One-dimensional Cramér type sharp asymptotics.* Let $\tau(x)$ denote the first passage time of the level $-x$ by $Z$

$$\tau(x) = \inf\{t \geq 0 : x + Z(t) < 0\} = \inf\{t \geq 0 : X(t) < 0\},$$

where $X$ is the translation of $Z$ by $x$, $X(t) = x + Z(t)$, and set $\psi(x) = P(\tau(x) < \infty)$. Let us assume the *Cramér assumption* that there exists a $\gamma > 0$, such that $\kappa(-\gamma) = 0$.

Under this assumption, Cramér has shown that when $Z$ is a Cramér–Lundberg process, then

(11) $$\lim_{x \to \infty} e^{\gamma x} \psi(x) = C,$$

where $C$ is given explicitly by [see, e.g., Feller (1971), Chapter XII.5]

(12) $$C = -\kappa'(0)/\kappa'_+(-\gamma),$$

which is strictly positive if and only if $\kappa'_+(-\gamma) > -\infty$. For a general Lévy process satisfying the Cramér assumption, Bertoin and Doney (1994) proved that (11) remains valid for some constant $C \geq 0$ that can be expressed in terms of the law of the ladder process and that is positive precisely if $E[|Z(1)|] < \infty$.

2.2. *The Arfwedson–Höglund theorem.* We recall now the asymptotics of the finite time ruin probability $\psi(x,t) = P(\tau(x) \leq t)$ for the *Cramér–Lundberg process*, first obtained by Arfwedson (1955) via the saddle-point method. Later, Höglund (1990) noted similar results for the probability of ruin after time $t$

(13) $$w(t,x) = P(t < \tau(x) < \infty) = \psi(x) - \psi(t,x).$$

Our formulation below, based on Höglund (1990), Corollary 2.3 and Asmussen (2000), uses the *exponential family* of measures $\{P^{(c)}\}$ defined for all $c \in \Theta = \{\theta : \kappa(\theta) < \infty\}$ by the Radon–Nikodym derivative $\Lambda^{(c)}$

$$(14) \quad \left.\frac{dP^{(c)}}{dP}\right|_{\mathcal{F}_t} = \Lambda^{(c)}(t) := \exp(c(X_t - X_0) - \kappa(c)t) = \exp(cZ_t - \kappa(c)t)$$

with $P^{(0)} = P$ and the corresponding shifted cumulant exponent given by

$$(15) \quad \kappa^{(c)}(\theta) := \kappa(\theta + c) - \kappa(c).$$

Further concepts, familiar from large deviations theory, that will be used are the convex conjugate $\kappa^*$ of the cumulant exponent $\kappa$, defined by

$$(16) \quad \kappa^*(v) = \sup_{\beta \in \mathbb{R}}[v\beta - \kappa(\beta)],$$

and the reparametrization of the exponential family by the corresponding set of means of $Z$. More precisely, it holds that to any $-v \in (\underline{v}, \overline{v})$ with

$$\underline{v} = \lim_{\theta \downarrow \underline{\theta}} \kappa'(\theta) \quad \text{and} \quad \overline{v} = \lim_{\theta \uparrow \overline{\theta}} \kappa'(\theta),$$

is associated a unique shift $\theta \in \Theta^o = (\underline{\theta}, \overline{\theta}) = (\underline{\theta}, \infty)$ such that

$$\kappa'(\theta_v) = -v.$$

Further, for any $0 < v < -\underline{v}$, the conjugate shift $\theta'_v \in \Theta$ is defined via

$$\kappa(\theta_v) = \kappa(\theta'_v), \qquad \theta_v < \theta'_v$$

and note that $v > 0$ implies $\theta_v < \theta'_v$. (Note that we parameterized by $-v$ since the possible means leading from $x > 0$ to $0$ are necessarily negative.)

REMARK 1. From now on, all quantities related to the shifted measure $P^{(c)}$ (like, e.g., the cumulant exponent, the adjustment coefficient, etc.) will be indicated by a superscript $(c)$ added to their $P$-counterparts. Observe that, in view of (15), the convex conjugate under $P^{(c)}$ and $P$ are related to each other by $\kappa^{*(c)}(v) = \kappa^*(v) + \kappa(c) - cv$. Similarly, it follows that the shift $\theta_v^{(c)}$ and its conjugate $\theta_v^{(c)\prime}$ are related to $\theta_v, \theta'_v$ by $\theta_v^{(c)} = \theta_v - c$ and $\theta_v^{(c)\prime} = \theta'_v - c$, respectively.

THEOREM 1. *Assume that either $\kappa'_+(0) < 0$ or that the Cramér assumption holds and write $\zeta = -\min\{\theta : \kappa(\theta) = 0\}$. If $0 < v < -\underline{v}$ and $x, t \to \infty$ such that $x/t = v$; it holds that*

$$(17) \quad \psi(x,t) \sim \begin{cases} Ce^{-x\zeta}, & \text{if } x/t < -\kappa'_+(-\zeta), \\ |D(v)|t^{-1/2}e^{-t\kappa^*(-v)}, & \text{if } x/t > -\kappa'_+(-\zeta), \end{cases}$$

$$(18) \quad w(x,t) \sim \begin{cases} |D(v)|t^{-1/2}e^{-t\kappa^*(-v)}, & \text{if } x/t < -\kappa'_+(-\zeta), \\ Ce^{-x\zeta}, & \text{if } x/t > -\kappa'_+(-\zeta), \end{cases}$$

with $C = 1$ if $\zeta = 0$ and $C = -\kappa'(0)/\kappa'_+(-\zeta)$ if $\zeta > 0$ and

$$(19) \qquad D(v) = c(v) \cdot \frac{1}{\sqrt{2\pi \kappa''(\theta_v)}} \qquad \text{with } c(v) = \frac{\theta'_v - \theta_v}{\theta_v \theta'_v},$$

where (17) is to be understood as $\lim_{x,t\to\infty, x=tv} e^{x\zeta} \psi(x,t) = 0$, if $\kappa'(-\zeta) = -\infty$ and $\zeta > 0$.

2.3. *The asymptotic limit laws of the process before/after ruin.* The following result shows existence of and identifies the limit laws of the Cramér–Lundberg process $X(t)$ conditioned on $t < \tau$ and on $t > \tau$. It will be used to establish Propositions 2 and 3.

THEOREM 2. *Assume that $\underline{\theta} < 0$ and $\kappa'(0) < 0$.*

(i) *If $0 < v < -\kappa'(0)$, then*

$$(20) \qquad \overline{\Psi}_v = \lim_{x,t\to\infty, x=tv} P(X(t) \in \cdot \,|\, \tau(x) > t),$$

*in the sense of weak convergence, where*

$$(21) \qquad \overline{\Psi}_v(dy) = c(v)^{-1}[e^{-\theta_v y} - e^{-\theta'_v y}]\mathbf{1}_{(0,\infty)}(y)\,dy$$

*with $c(v)$ the function appearing in (19).*

(ii) *If $-\kappa'(0) < v < -\underline{v}$, then*

$$(22) \qquad \Psi_v = \lim_{x,t\to\infty, x=tv} P(X(t) \in \cdot \,|\, \tau(x) < t)$$

*in the sense of weak convergence, where*

$$(23) \qquad \Psi_v(dy) = |c(v)|^{-1}[e^{-\theta'_v y}\mathbf{1}_{(0,\infty)}(y) + e^{-\theta_v y}\mathbf{1}_{(-\infty,0)}(y)]\,dy.$$

PROOF. (i) First, we verify that the measure $\overline{\Psi}_v$ is a probability measure. Indeed, it is not hard to verify that $\overline{\Psi}_v$ is a measure (since $\theta'_v > \theta_v > 0$) that integrates to one. Further, it is easily checked that the mgf $M^v$ of $\overline{\Psi}_v$ is given by

$$M^v(c) = \frac{\theta_v}{\theta_v - c} \cdot \frac{\theta'_v}{\theta'_v - c} \qquad \text{for } c < \theta_v.$$

In view of the continuity theorem, the weak convergence in (20) follows once we show that the mgfs $M_{x,t}(c)$ of the measures $P(X(t) \in \cdot \,|\, \tau(x) > t)$ converge pointwise to $M^v(c)$ as $x, t \to \infty, x/t = v$, for $c$ in some neighborhood of the origin.

Since $\underline{\theta} < 0$, it holds that for all $c$ in some neighborhood of the origin $M_{x,t}(c)$ is finite and $E[e^{cX(t)}\mathbf{1}_{\{t<\tau(x)\}}] = e^{cx+\kappa(c)t}P^{(c)}(t<\tau(x))$. It also holds that $P^{(c)}(t<\tau(x)) = P^{(c)}(t<\tau(x)<\infty)$ [since $\kappa^{(c)\prime}(0^+) < 0$]. Therefore,

$$M_{x,t}(c) = E[e^{cX(t)}|t<\tau(x)] = \frac{E[e^{cX(t)}\mathbf{1}_{\{t<\tau(x)<\infty\}}]}{P(t<\tau(x)<\infty)}$$

$$= \frac{e^{cvt+\kappa(c)t}P^{(c)}(t<\tau(x)<\infty)}{P(t<\tau(x)<\infty)}.$$

Invoking Theorem 1 and (15) for the numerator and denominator, it follows by taking the limit of $x,t \to \infty$, $x = tv$ that

$$\lim_{x,t\to\infty, x=tv} M_{x,t}(c) = D^{(c)}(v)/D(v) = \frac{\theta_v^{(c)\prime} - \theta_v^{(c)}}{|\theta_v^{(c)}\theta_v^{(c)\prime}|}\frac{|\theta_v^\prime\theta_v|}{\theta_v^\prime - \theta_v}.$$

In view of Remark 1, the latter is equal to $M^v(c)$.

The proof of (ii) is similar and omitted. $\square$

2.4. *Law of large numbers for the ruin time.* We include now for reference a result concerning the behavior of the time of ruin of *a general Lévy process $Z$* for large initial reserves.

LEMMA 1. *Suppose that $E[|Z(1)|] < \infty$ and $E[Z(1)] \leq 0$. Then, as $x \to \infty$:*

(i) $\tau(x)/x \to -E[Z(1)]^{-1}$ *P-a.s. and*
(ii) $E[\tau(x)]/x \to -E[Z(1)]^{-1}$.

PROOF. (ii) If $E[Z(1)] = 0$, then the Lévy process $Z$ oscillates and the identity follows since then $E[\tau(x)] = +\infty$ for every $x$ (see, e.g., Bertoin (1996), Chapter VI, Proposition 17(iii)). Suppose now that $-\infty < E[Z(1)] < 0$ (so that $Z$ drifts to $-\infty$) and first exclude the case that $Z$ is a compound Poisson process. Denoting by $L^{-1}(t) = \inf\{u \geq 0 : L(u) > t\}$, the inverse of the local time $L$ of $Z$ and $T(x) = \inf\{t \geq 0 : H(t) > x\}$ the first passage time of the ladder height process $H(t) = Z(L^{-1}(t))$ of $Z$ it is easily verified that $\tau(x) = L^{-1}(T(x))$. The pair $(L^{-1}, H)$ forms a two-dimensional Lévy process and we denote its bivariate Laplace exponent by $\widehat{\kappa}$. The Laplace transform of $E[\tau(x)]$ can then be expressed as follows:

$$\int_0^\infty e^{-\lambda x} E[\tau(x)]\,dx = \frac{\partial_1\widehat{\kappa}(0, 0^+)}{\lambda\widehat{\kappa}(0, \lambda)},$$

where $\partial_i$ denotes the partial derivative with respect to $i$th variable (see, e.g., Bertoin (1996), Chapter VI, Proposition 17). As $\widehat{\kappa}(0,0) = 0$ and $\widehat{\kappa}(0,\cdot)$ is right-differentiable in zero, it follows in view of a Tauberian theorem that

$$E[\tau(x)] \sim x\,\partial_1\widehat{\kappa}(0, 0^+)/\partial_2\widehat{\kappa}(0, 0^+) \qquad \text{as } x \to \infty.$$

The strong law of large numbers implies that the product $H(t)/L^{-1}(t) = [Z(L^{-1}(t))/t] \times [t/L^{-1}(t)]$ converges to

$$E[Z(1)] = E[Z(L^{-1}(1))]E[L^{-1}(1)]^{-1} \qquad (24)$$

(the corresponding result for random walks is known as the famous *Wald identity*). Since $\partial_1 \widehat{\kappa}(0^+, 0) = E[L^{-1}(1)]$ and $\partial_2 \widehat{\kappa}(0, 0^+) = E[Z(L^{-1}(1))]$, the claim follows. The case of a compound Poisson process follows by adding a small drift.

(i) The strong law of large numbers implies that, $P$-a.s.,

$$\tau(x)/x = L^{-1}(T(x))/T(x) \cdot T(x)/x \to E[L^{-1}(1)]/E[H(1)] = E[Z(1)],$$

as $x \to \infty$, where we used the Wald identity (24). □

**3. The exact ultimate ruin probability for the degenerate 2-d process.** In this section, we consider the probability that a Lévy process $S$ starting at $0$ ever upcrosses a piecewise linear barrier $b$. To be specific, we consider the first passage time $\tau_b$ of $S$ over $b$, as in (8), where $b$ is given by $b = b_{\min}$ or $b = b_{\max}$ with

$$b_{\min}(t) = \min_{i=1,2}\{x_i + p_i t\}, \qquad b_{\max}(t) = \max_{i=1,2}\{x_i + p_i t\}$$

where $x_2 > x_1$ and $p_1 > p_2$. As noted in the Introduction, $P(\tau_b < \infty)$ with $b = b_{\min}$ (resp. $b = b_{\max}$) exactly coincides with the ruin probability $\psi_{\text{or}}(x_1, x_2)$ [resp. $\psi_{\text{sim}}(x_1, x_2)$] of the process $(X_1, X_2)$ with

$$X_i(t) := x_i + p_i t - S(t), \qquad i = 1, 2. \qquad (25)$$

Denoting by

$$T = T(x_1, x_2) = \frac{x_2 - x_1}{p_1 - p_2} \qquad (26)$$

the time at which the lines $t \mapsto x_1 + p_1 t$ and $t \mapsto x_2 + p_2 t$ cross, we see that, for example, for $S$ to never cross $b_{\min}$, it is required to stay below the barrier $x_1 + p_1 t$ between the times $0$ and $T$ and subsequently to stay below the barrier $x_2 + p_2 t$ after time $T$. Since $S$ is Markovian and $x_1 + p_1 T = x_2 + p_2 T$, conditioning at time $T$ yields

$$\overline{\psi}_{\text{or}}(x_1, x_2) = \int_0^\infty f_1(ds, T|x_1) F_2(s),$$

where $\overline{\psi}_{\text{or}} = 1 - \psi_{\text{or}}$, $F_i(s) = P(s + S(t) \leq x_i + p_i(t + T) \ \forall t > 0)$ and

$$f_i(ds, T|x) := P(S(t) \leq x + p_i t, \ \forall t \in [0, T], S(T) \in ds) \qquad (27)$$

is the density of $S(T)$ of the paths at time $T$ that "survived" the upper barrier $x + p_i t$. Reformulating this result in terms of the two coordinates $X_i$ of $X = (X_1, X_2)$ and the coordinate-wise densities of the surviving paths

$$\overline{\psi}_i(dz, T|x_i) = P(\tau_i(x_i) > T, X_i(T) \in dz) \qquad (28)$$

we arrive thus at the following result, which relates the ruin probabilities of the two dimensional process $X$ to those of its coordinates $X_1, X_2$.

PROPOSITION 1. *Let $X = (X_1, X_2)$ be the two-dimensional Lévy process with $X_i$ given in (25) and suppose that $x_2 > x_1$ and $p_2 < p_1$.*

(a) *The ruin probabilities $\psi_{\text{and}}, \psi_{\text{sim}}$ and $\psi_{\text{or}}$ are given by*

$$\psi_{\text{sim}}(x_1, x_2) = P(\tau_2(x_2) \leq T) + P\Big(\tau_2(x_2) > T, \inf_{s>T} X_1(s) < 0\Big), \tag{29}$$

$$\psi_{\text{or}}(x_1, x_2) = P(\tau_1(x_1) \leq T) + P\Big(\tau_1(x_1) > T, \inf_{s>T} X_2(s) < 0\Big), \tag{30}$$

$$\psi_{\text{and}}(x_1, x_2) = P(T < \tau_1(x_1) < \infty) + P(\tau_1(x_1) \leq T, \tau_2(x_2) < \infty). \tag{31}$$

(b) *The survival probabilities $\overline{\psi}_{\text{or}} = 1 - \psi_{\text{or}}$ and $\overline{\psi}_{\text{sim}} = 1 - \psi_{\text{sim}}$ are given by*

$$\overline{\psi}_{\text{sim}}(x_1, x_2) = P\Big(\tau_2(x_2) > T, \inf_{s>T} X_1(s) \geq 0\Big) = \int_0^\infty \overline{\psi}_2(dz, T|x_2)\overline{\psi}_1(z),$$

$$\overline{\psi}_{\text{or}}(x_1, x_2) = P\Big(\tau_1(x_1) > T, \inf_{s>T} X_2(s) \geq 0\Big) = \int_0^\infty \overline{\psi}_1(dz, T|x_1)\overline{\psi}_2(z),$$

*where $T$ is given in (26), $\overline{\psi}_i(dz, T|x_i)$ in (28) and $\overline{\psi}_i(z) = P(\tau_i(z) = \infty)$ are perpetual one-dimensional survival probabilities.*

PROOF. By definition of $\psi_{\text{sim}}$, it holds that

$$\overline{\psi}_{\text{sim}}(x_1, x_2) = P(\max\{X_1(t), X_2(t)\} \geq 0 \text{ for all } t \geq 0).$$

Next, we note that, if $x_2 > x_1$, it holds that the maximum

$$\max\{X_1(t), X_2(t)\} = \max\{x_1 - x_2 + (p_1 - p_2)t, 0\} + X_2(t)$$

is equal to $X_2(t)$ for $t \leq T$ and to $X_1(t)$ for $t > T$, where $T$ was defined in (26). Applying subsequently the Markov property of $X_1$ at time $T$ shows that

$$\begin{aligned}\overline{\psi}_{\text{sim}}(x_1, x_2) &= P(X_2(t) \geq 0 \text{ for } t \leq T, X_1(t) \geq 0 \text{ for } t \geq T) \\ &= \int_0^\infty P(X_2(T) \in dz, \tau(x_2) > T)P(\tau_1(z) = \infty).\end{aligned} \tag{32}$$

The identity in (29) follows by taking the complement of (32). The proof of $\psi_{\text{or}}$ is similar and omitted. Finally, write

$$\psi_{\text{and}}(x_1, x_2) = P(T < \tau_1(x_1) < \infty, \tau_2(x_2) < \infty) + P(\tau_1(x_1) \leq T, \tau_2(x_2) < \infty).$$

Equation (31) follows then by checking that $\{\tau_1(x_1) > T\}$ and $\{\inf_{s>T} X_1(s) < 0\}$ respectively imply that $\{\tau_2(x_2) > T\}$ and $\{\inf_{s>T} X_2(s) < 0\}$. □

In the special case that $S$ is a compound Poisson process with exponential claims $\sigma_i$ with parameter $\mu$, we have exponential ultimate ruin probabilities

$$\psi_i(x) = C_i e^{-\gamma_i x}, \tag{33}$$

where $C_i = \frac{\lambda}{\mu p_i}$ and $\gamma_i = \mu - \lambda/p_i$. Similarly, if $S$ is a spectrally negative Lévy process, the Markov property and the absence of positive jumps imply the multiplicativity property $P(\tau_i(x+y) < \infty) = P(\tau_i(x) < \infty)P(\tau_i(y) < \infty)$. Thus, $P(\tau_i(x) < \infty)$ must be an exponential function (33) and the constant $C_i$ equals 1. For these two cases, equations (29)–(31) can be developed further by employing the technique of change of measure.

3.1. *The case of exponential ultimate ruin probabilities.* Consider now the relation (30) in the case of exponential ultimate ruin probabilities that is when $\psi_i(x) = C_i e^{-\gamma_i x}$ for $C_i, \gamma_i > 0$ ($i = 1, 2$). Note that

$$\psi_{\text{or}}(x_1, x_2) = P(\tau_1(x_1) \leq T) + C_2 E[e^{-\gamma_2 X_1(T)} \mathbf{1}_{\{\tau_1(x_1) > T\}}], \tag{34}$$

where $C_2 = 1$ in the case that $S$ is a spectrally negative Lévy process. Let $\kappa_i$ be the cumulant exponent of $X_i(t) - x_i$. By a change of measure (14) and using that $-\gamma_2 x_1 + \kappa_1(-\gamma_2)T = -\gamma_2 x_2$, we find that the second term in (34) is equal to

$$C_2 e^{-\gamma_2 x_1 + \kappa_1(-\gamma_2)T} E[\Lambda^{(-\gamma_2)}(T) \mathbf{1}_{\{\tau_1(x_1) > T\}}] = C_2 e^{-\gamma_2 x_2} P^{(-\gamma_2)}(\tau_1(x_1) > T).$$

The probabilities $\psi_{\text{sim}}$ can be treated using similar arguments, and $\psi_{\text{and}}$ is obtained from the "complementarity equation" (7). In conclusion, the original two-dimensional ruin problems $\psi_{\text{or}}/\psi_{\text{sim}}/\psi_{\text{and}}$ are reduced to one-dimensional finite time ruin problems $\psi_i^{(c)}(x,t) = P^{(c)}(\tau_i(x) \leq t)$, as follows.

COROLLARY 1. *Suppose $S$ is a spectrally negative Lévy process, or a compound Poisson process with exponential (positive) jumps and let $X_i$ defined by (25). If $x_2 > x_1$, it holds that*

$$\psi_{\text{sim}}(x_1, x_2) = \psi_2(x_2, T) + \psi_1(x_1)\overline{\psi}_2^{(-\gamma_1)}(x_2, T),$$

$$\psi_{\text{or}}(x_1, x_2) = \psi_1(x_1, T) + \psi_2(x_2)\overline{\psi}_1^{(-\gamma_2)}(x_1, T),$$

$$\psi_{\text{and}}(x_1, x_2) = w_1(x_1, T) + \psi_2(x_2)\psi_1^{(-\gamma_2)}(x_1, T),$$

*where $\overline{\psi}_i^{(c)}(x,t) = 1 - \psi_i^{(c)}(x,t)$ and*

$$w_1(x,t) = P(t < \tau_1(x) < \infty) = \psi_1(x) - \psi_1(x,t). \tag{35}$$

PROOF. Let us establish the last statement. By inserting the expression for $\psi_{\text{or}}$ in (7), it follows that

$$\begin{aligned}\psi_{\text{and}}(x_1, x_2) &= \psi_1(x_1) + \psi_2(x_2) - \psi_{\text{or}}(x_1, x_2) \\ &= \psi_1(x_1) - \psi_1(x_1, T) + \psi_2(x_2)(1 - \overline{\psi}_1^{(-\gamma_2)}(x_1, T)). \quad \square\end{aligned}$$

This decomposition result (and its generalization) will provide the key for obtaining the two terms asymptotic expansions in Propositions 2, 3 below.

**4. General two-dimensional Cramér asymptotics.** We consider now the asymptotics of the ruin probabilities $\psi_{\text{sim}}, \psi_{\text{or}}$ and $\psi_{\text{and}}$ when the initial reserves tend to infinity along a ray, for a general two-dimensional Lévy process $X = (X_1, X_2)$ starting from $x = (x_1, x_2)$. To avoid degeneracies, we exclude throughout the cases that $X_1$ or $X_2$ have monotone paths, or that the ratio $[X_1 - x_1]/[X_2 - x_2]$ is constant. The law of the process $X$ is determined by its joint cumulant exponent $\kappa(\theta_1, \theta_2) = \log E[e^{\theta_1(X_1(1)-x_1) + \theta_2(X_2(1)-x_2)}]$ which is well defined on its domain $\Xi = \{\theta \in \mathbb{R}^2 : \kappa(\theta) < \infty\}$, whose interior is denoted by $\Xi^o$. For every $\theta \in \Xi^o$, the gradient $\nabla \kappa(\theta) = (\partial_1 \kappa(\theta), \partial_2 \kappa(\theta))$ is well defined. Other subsets of $\Xi$ playing a role in our setting are: the *Cramér set* $\mathcal{C}$, its interior $\mathcal{C}^o$ and its boundary $\partial \mathcal{C} := \mathcal{C} \setminus \mathcal{C}^o$ where

$$\mathcal{C} = \{(\theta_1, \theta_2) \in \Xi : \kappa(\theta_1, \theta_2) \leq 0\}.$$

In view of the convexity of $\kappa$, it follows that the set $\mathcal{C}$ is convex and that for fixed $\theta' \in \partial \mathcal{C} \cap \Xi^o$, it holds that

(36) $$[\theta - \theta'] \cdot \nabla \kappa(\theta') \leq 0 \qquad \text{for all } \theta \in \mathcal{C},$$

where $\cdot$ denotes the inner-product.

Associated to any $c \in \Xi$ is a measure $P^{(c)}$ defined as a twist of $P$ by the martingale $\exp(c_1(X_1(t) - x_1) + c_2(X_2(t) - x_2) - \kappa(c_1, c_2)t)$.

We assume throughout that besides the origin, the Cramér set intersects the axes in two more points: $\boldsymbol{\gamma}^{(1)} = (-\gamma_1, 0), \boldsymbol{\gamma}^{(2)} = (0, -\gamma_2) \in \partial \mathcal{C}, \gamma_1, \gamma_2 > 0$, so that

(37) $$\kappa(-\gamma_1, 0) = \kappa(0, -\gamma_2) = 0.$$

We shall also assume that

(38) $$(-\gamma_1, 0), (0, -\gamma_2) \in \Xi^o.$$

When (37), (38) hold, we will say that the Cramér assumptions hold true for $X = (X_1, X_2)$.

EXAMPLE 1. If $S(t)$ is a Lévy process and $X_i(t) = x_i + p_i t - S(t)$ (with the cumulant generating functions $\kappa_i$), then the joint cumulant generating function $\kappa$ of $(X_1, X_2)$ is related to $\kappa_1$ and $\kappa_2$ by

$$\kappa(\theta_1, \theta_2) = \kappa_1(\theta_1 + \theta_2) - \theta_2(p_1 - p_2) = \kappa_2(\theta_1 + \theta_2) + \theta_1(p_1 - p_2)$$
$$= p \cdot \theta + \kappa_S(\theta_1 + \theta_2),$$

where $\kappa_S$ is the cumulant exponent of $S$. It is easy to check that the degenerate two-dimensional Lévy process $X = (X_1, X_2)$ satisfies the Cramér-conditions iff its coordinates do, that is, if there exist constants $\gamma_i > 0$, in the interior of the domains of the cumulant exponents $\kappa_i$ of $X_i(t) - x_i$ ($i = 1, 2$), such that

$$\kappa_i(-\gamma_i) = 0. \tag{39}$$

The following result yields the asymptotics of $\psi_{\text{or}}$ and an order estimate of $\psi_{\text{and}}$ for general two-dimensional Lévy processes.

THEOREM 3. *Suppose that the Cramér assumptions (37) and (38) hold, and let $a > 0$. Then as $K \to \infty$,*

$$\psi_{\text{or}}(aK, K) \sim C_2 e^{-\gamma_2 K} + C_1 e^{-\gamma_1 a K}, \tag{40}$$

$$\psi_{\text{and}}(aK, K) = o(C_2 e^{-\gamma_2 K} + C_1 e^{-\gamma_1 a K}), \tag{41}$$

*where $C_i > 0, i = 1, 2$, are the asymptotic constants corresponding to $X_i$.*

The proof of Theorem 3 is based on the following estimates.

LEMMA 2. *The following hold true:*

(i) $\max\{\psi_1(x_1), \psi_2(x_2)\} \leq \psi_{\text{or}}(x_1, x_2) \leq \psi_1(x_1) + \psi_2(x_2)$;
(ii) $\psi_{\text{or}}(x_1, x_2) = \psi_{1 \leq 2}(x_1, x_2) + \psi_{2 \leq 1}(x_1, x_2) - \psi_{1=2}(x_1, x_2)$, *where*

$$\psi_{i \leq j}(x_1, x_2) := P(\tau_i(x_i) \leq \tau_j(x_j), \tau_i(x_i) < \infty)$$

*and*

$$\psi_{1=2}(x_1, x_2) := P(\tau_1(x_1) = \tau_2(x_2) < \infty).$$

PROOF. The estimates follow in view of the observations that

$$\{\tau_i(x_i) < \infty\} \subset \{\tau_{\text{or}}(x_1, x_2) < \infty\} \subset \bigcup_{i=1}^{2} \{\tau_i(x_i) < \infty\} \quad \text{for } i = 1, 2,$$

and

$$\{\tau_{\text{or}}(x_1, x_2) < \infty\} = A_1 \cup A_2 \setminus [A_1 \cap A_2],$$

where $A_i = \{\tau_i(x_i) < \infty, \tau_i(x_i) \leq \tau_{3-i}(x_{3-i})\}$. □

LEMMA 3. *Suppose that (37) and (38) hold, and write $\gamma_{a,\beta} := \beta a\gamma_1 + (1-\beta)\gamma_2$.*

(i) *If $a\gamma_1 = \gamma_2$ it holds that, as $K \to \infty$,*

(42) $\qquad \psi_{1\leq 2}(aK, K) \sim C_1 e^{-\gamma_1 aK}, \qquad \psi_{2\leq 1}(aK, K) \sim C_2 e^{-\gamma_2 K}.$

(ii) *For $a > 0$ and any $\beta \in (0,1)$, $\psi_{\text{sim}}(aK, K) = \text{o}(e^{-\gamma_{a,\beta}K})$ $(K \to \infty)$.*

This lemma (established below) implies immediately Theorem 3.

PROOF OF THEOREM 3. First note that, in view of the Cramér–Lundberg asymptotics (11) and equation (7), the asymptotics in (40) imply the estimate in (41). The rest of the proof is therefore devoted to establishing (40).

In view of (11) and Lemma 2(i), it follows that, if $\gamma_1 a > \gamma_2$ [resp. $\gamma_1 a < \gamma_2$], the lower bound and upper bound in Lemma 2(i) are of the same order of magnitude, $C_2 e^{-\gamma_2 K}$ [resp. $C_1 e^{-\gamma_1 aK}$], as $K \to \infty$. Thus, (40) is valid if $\gamma_1 a \neq \gamma_2$.

Next we turn to the case $\gamma_1 a = \gamma_2$. Since $\psi_{1=2}$ is dominated by $\psi_{\text{sim}}$ and $\gamma_{a,\beta} = \gamma_2 = a\gamma_1$ if $a\gamma_1 = \gamma_2$, it follows, by invoking Lemma 3(ii), that $\psi_{1=2}(aK, K) = \text{o}(e^{-\gamma_2 K}) = \text{o}(e^{-\gamma_1 aK})$ as $K \to \infty$. In view of Lemma 2(ii) and Lemma 3(i), it therefore follows that (40) is also valid if $a\gamma_1 = \gamma_2$. $\square$

PROOF OF LEMMA 3. (i) The asymptotics of $\psi_{1\leq 2}$ follow once we have shown that as $K \to \infty$ it holds that

(43) $\quad e^{\gamma_1 aK}\psi_{1\leq 2}(aK, K) = E^{(-\gamma_1,0)}(e^{-\gamma_1 X_1(\tau_1)}\mathbf{1}_{\{\tau_1 \leq \tau_2, \tau_1 < \infty\}}) \to C_1,$

where $\tau_1 = \tau_1(aK)$ and $\tau_2 = \tau_2(K)$. To prove this claim, we compare the asymptotic behavior of $\tau_1$ and $\tau_2$ as $K \to \infty$, adapting the argument developed in Glasserman and Wang (1997) (Proposition 2) for random walk. If $E^{(-\gamma_1,0)}[X_2(1)] > x_2$, then $P^{(-\gamma_1,0)}(\tau_2 = \infty) \to 1$ as $K \to \infty$ and, invoking (11), the claim (43) follows. If $E^{(-\gamma_1,0)}[X_2(1)] \leq x_2$, it follows in view of Lemma 1(i) that as $K \to \infty$, $P^{(-\gamma_1,0)}$-a.s.,

(44) $\qquad \dfrac{\tau_1(aK)}{\tau_2(K)} = a\dfrac{\tau_1(aK)}{aK}\dfrac{K}{\tau_2(K)} \to \dfrac{a\partial_2\kappa(-\gamma_1,0)}{\partial_1\kappa(-\gamma_1,0)},$

where we used that $\partial_i\kappa(\theta) = E^{(\theta)}[X_i(1) - x_i]$ for $\theta \in \Xi^o$, $i=1,2$. Applying (36) with $\theta = (0,-\gamma_2)$ and $\theta' = (-\gamma_1, 0)$, we see that the right-hand side of (44) is bounded above by $a\gamma_1/\gamma_2$, which is equal to one if $\gamma_2 = a\gamma_1$. Therefore, $\tau_2(K)$ dominates $\tau_1(aK)$ for all $K$ large enough and (43) follows as a consequence of the Cramér–Lundberg asymptotics (11). The asymptotics of $\psi_{2\leq 1}$ can be treated similarly.

(ii) Choose $\beta \in (0,1)$ and write $\gamma(\beta) = \beta(\gamma_1, 0) + (1-\beta)(0, \gamma_2) = (\gamma_1(\beta), \gamma_2(\beta))$. The key step is to verify that the segment $\{\gamma(\beta), \beta \in (0,1)\}$ is not part of the boundary $\partial C$. Indeed, since the function $\beta \mapsto f(\beta)$ with $f(\beta) = \exp(\kappa(-\gamma(\beta))) = E[\exp\{\beta\gamma_1(X_1(1)-x_1)+(1-\beta)\gamma_2(X_2(1)-x_2)\}]$ is strictly convex with $f(0) = f(1) = 1$, it follows that $f(\beta) < 1$ for $\beta \in (0,1)$. (The strict convexity is a consequence of the facts that $f''(\beta) = E[(\gamma_1(X_1(1) - x_1) - \gamma_2(X_2(1) - x_2))^2 e^{\beta\gamma_1(X_1(1)-x_1)+(1-\beta)\gamma_2(X_2(1)-x_2)}] > 0$ and we excluded the case $[X_1 - x_1]/[X_2 - x_2] = \text{const.}$)

Therefore, there exists a $-\gamma^* = -(\gamma_1^*, \gamma_2^*) \in \mathcal{C}^o$ such that $\gamma_i^* > \gamma_i(\beta)$ ($i = 1, 2$). By changing the measure, we see that $\psi_{\text{sim}}(aK, K)$ is equal to

$$e^{-(\gamma_1^* a + \gamma_2^*)K} E^{(-\gamma^*)}[e^{\gamma_2^* X_2(\tau_{\text{sim}}) + \gamma_1^* X_1(\tau_{\text{sim}}) + \kappa(-\gamma_1^*, -\gamma_2^*)\tau_{\text{sim}}} \mathbf{1}_{\{\tau_{\text{sim}} < \infty\}}],$$

where $\tau_{\text{sim}} = \tau_{\text{sim}}(aK, K)$. Since $X_i(\tau_{\text{sim}}) \leq 0$ and $\kappa(-\gamma_1^*, -\gamma_2^*) \leq 0$, this expectation is bounded above by 1, and as $a\gamma_1^* + \gamma_2^* > \gamma_{a,\beta}$, it thus follows that $\psi_{\text{sim}}(aK, K) = \text{o}(e^{-\gamma_{a,\beta}K})$ as $K \to \infty$. □

The following result concerns the asymptotics in the upper cone $\{x_1 \leq x_2\}$, in the case of the Sparre–Andersen model.

THEOREM 4. *Let $S$ be a compound renewal process as in (2) and let $a < 1$. Assume there exist $\gamma_i > 0$ and an $\epsilon > 0$ such that $E[e^{-\gamma_i p_i \zeta}]E[e^{\gamma_i \sigma}] = 1$ and $E[e^{(\gamma_i + \epsilon)\sigma}] < \infty$. Then it holds that, as $K \to \infty$,*

(45) $$\psi_{\text{or}}(aK, K) \sim C_2 e^{-\gamma_2 K} + C_1 e^{-\gamma_1 aK},$$

(46) $$\psi_{\text{and}}(aK, K) = \text{o}(C_2 e^{-\gamma_2 K} + C_1 e^{-\gamma_1 aK}).$$

PROOF. In view of the key observation that the ruin probabilities $\psi_{\text{or}}/\psi_{\text{and}}$ and $\psi_{\text{sim}}$ do not change if we replace $X = (X_1, X_2)$ by a two-dimensional compound Poisson process with unit jump rate and jump sizes distributed as $(\sigma_n - p_1\zeta_n, \sigma_n - p_2\zeta_n)$, the statement follows by invoking Theorem 3. □

In Section 5 below, we will sharpen Theorem 3, in the degenerate case. Before that, we introduce a partition of the quadrant in cones, which turn out to describe the different asymptotic regimes of the ruin probabilities as the initial reserves $(x_1, x_2)$ tend to infinity along a ray.

4.1. *The asymptotic cones.* We introduce now two cones $\mathcal{D}_i$ ($i = 1, 2$) within the quadrant $\mathbb{R}_+^2 = (0, \infty)^2$, situated between the $x_i$ axis and the directions of the expected drift evaluated at the adjustment tilts $\mathbf{v}^{(i)} := \nabla\kappa(\boldsymbol{\gamma}^{(i)}), i = 1, 2$. Heuristically, these cones can be described as the "asymptotic boundary cones" of the "sim" ruin, that is, the cones where asymptotically this event happens dominantly by straight paths running to the half-lines $\{x : x_i = 0, x_{3-i} \leq 0, i = 1, 2\}$.

Letting $s_i$ denote the "slopes $\frac{dx_1}{dx_2}$" of the vector $\nabla\kappa(\boldsymbol{\gamma}^{(i)})$ ($i = 1, 2$), the cones are given by

$$\mathcal{D}_1 = \{(x_1, x_2) \in \mathbb{R}_+^2 : x_1 > x_2 s_1\},$$

(47)
$$\mathcal{D}_2 = \{(x_1, x_2) \in \mathbb{R}_+^2 : x_1 < x_2 s_2\}.$$

Note that in the degenerate case, $s_i$ become

(48)
$$s_1 = \frac{\kappa_1'(-\gamma_1)}{\kappa_2'(-\gamma_1)}, \qquad s_2 = \left(\frac{\kappa_1'(-\gamma_2)}{\kappa_2'(-\gamma_2)}\right)_+$$

and we show in Lemma 4 that $s_2 < s_1$, implying that the cones $\mathcal{D}_1$ and $\mathcal{D}_2$ are disjoint in this case. Let now

$$\mathcal{D}_0 = \mathbb{R}_+^2 \setminus [\overline{\mathcal{D}}_1 \cup \overline{\mathcal{D}}_2]$$

denote the open cone lying between $\mathcal{D}_1$ and $\mathcal{D}_2$ (where $\overline{\mathcal{D}}_i$ is the closure of $\mathcal{D}_i$).

We will show in Theorem 5 that within $\mathcal{D}_i$, $\psi_{\text{sim}}(x_1, x_2)$ is asymptotically equivalent to $\psi_i(x_i), i=1,2$, respectively, and that a different regime holds within $\mathcal{D}_0$, which is characterized by "radial dependence" on the slope $a = a(x) := \frac{x_1}{x_2}$.

We will also show in Theorem 6 that in the case when $\mathcal{D}_2$ is void [which is characterized by $\partial_1 \kappa_1(\boldsymbol{\gamma}^{(2)}) > 0$], a special type asymptotic regime holds for $\psi_{\text{and}}$, within a new "secondary cone" $\widehat{\mathcal{D}}_2$ situated between the $x_2$ axis and the direction $\mathbf{v}^{(3)} := \nabla \kappa(\boldsymbol{\gamma}^{(3)})$, where $\boldsymbol{\gamma}^{(3)}$ is defined as the leftmost intersection of the Cramér set with the line $\theta_2 = -\gamma_2$. In the degenerate case, we have

$$\widehat{\mathcal{D}}_2 = \{(x_1, x_2) \in \mathbb{R}_+^2 : x_1 < x_2 s_3\}$$

where

$$s_3 = \kappa_1'(-\gamma_3)/\kappa_2'(-\gamma_3)$$

with $\gamma_3$ the largest root of $\kappa_1(-s) = \kappa_1(-\gamma_2)$. As stated in Lemma 4 below, if $\kappa_1'(-\gamma_2) > 0$, it follows that $\widehat{\mathcal{D}}_2 \neq \varnothing = \mathcal{D}_2$ and otherwise the cones $\mathcal{D}_2$ and $\widehat{\mathcal{D}}_2$ coincide. This partition of the positive quadrant into cones

$$\mathcal{D}_1, \widehat{\mathcal{D}}_2 \text{ and } \widehat{\mathcal{D}}_0 := \mathbb{R}_+^2 \setminus [\overline{\mathcal{D}}_1 \cup \overline{\widehat{\mathcal{D}}}_2]$$

is illustrated in Figure 2.

**5. Sharp asymptotics for degenerate risk-processes.** We restrict now ourselves to a two-dimensional Lévy process $(X_1, X_2)$ with $X_i(t) = x_i + p_i t - S(t)$, $i=1,2$, where $S$ is a compound Poisson process with positive jumps. Throughout this section, we assume that there exist $\gamma_i > 0$, $\quad i = 1, 2$, such that $\kappa_i(-\gamma_i) = 0$, where $\kappa_i$ is the cumulant exponent of $X_i - x_i$. Note that whenever $\gamma_1$ exists, $\gamma_2$ exists as well, and $\gamma_1 > \gamma_2$ (by the convexity of $\kappa_i$).

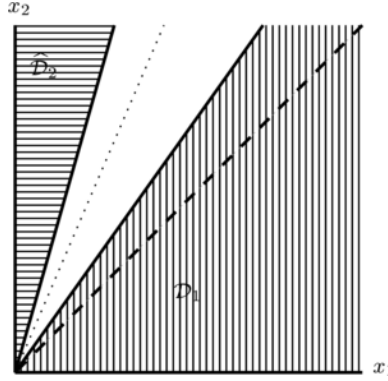

FIG. 2. *Pictured are the positive quadrant, divided in the three cones $\mathcal{D}_1$ (shaded with vertical lines), $\mathcal{D}_2$ or $\widehat{\mathcal{D}}_2$ (shaded with horizontal lines) and $\mathcal{D}_0$ or $\widehat{\mathcal{D}}_0$ (white), as well as the lines $x_1 = x_2$ (dashed) and $\gamma_1 x_1 = \gamma_2 x_2$ (dotted).*

5.1. *A characterization of the asymptotic cones.* We start with gathering some properties of the asymptotic cones, in the degenerate setting.

LEMMA 4. *The following hold true:*

(i) *The cones $\mathcal{D}_i, i = 0, 1, 2$ are disjoint and $\mathcal{D}_0, \mathcal{D}_1 \neq \varnothing$.*
(ii) *$\widehat{\mathcal{D}}_2 = \mathcal{D}_2 \neq \varnothing$ iff $\kappa_1'(-\gamma_2) < 0$ and $\mathcal{D}_2 = \varnothing \neq \widehat{\mathcal{D}}_2$ iff $\kappa_1'(-\gamma_2) > 0$.*
(iii) *$\mathcal{D}_1 \subset \mathcal{U} := \{(x_1, x_2) \in \mathbb{R}_+^2 : x_2 \gamma_2 < x_1 \gamma_1\}$ and $\mathcal{D}_2 \subset \mathbb{R}_+^2 \setminus \mathcal{U}$.*

PROOF. Writing

$$\frac{\kappa_1'(s)}{\kappa_2'(s)} = \frac{\kappa_2'(s) + p_1 - p_2}{\kappa_2'(s)} = 1 + \frac{p_1 - p_2}{\kappa_2'(s)},$$

it follows that $s_1 < 1$, since $\kappa_2'(-\gamma_1) < 0$, and that $s_2 < s_1$, since $\gamma_1 > \gamma_2$ and, by the strict convexity of $\kappa_2$, $\kappa_2'$ is strictly increasing on its domain. Next, in view of the definitions of $s_2$ and $s_3$, it follows that $s_3 = 0$ [resp. $s_2 = 0$] iff $\kappa_1'(-\gamma_2) = 0$ [resp. $\kappa_1'(-\gamma_2) \geq 0$]. Subsequently, we note that on the ray $x_1/x_2 = \gamma_2/\gamma_1$ it holds that

$$\frac{x_2}{T(x_1, x_2)} = \frac{p_1 - p_2}{1 - \gamma_2/\gamma_1} = \frac{\kappa_2(-\gamma_1) - \kappa_1(-\gamma_1)}{\gamma_1 - \gamma_2} = \frac{\kappa_2(-\gamma_1) - \kappa_2(-\gamma_2)}{\gamma_1 - \gamma_2}.$$

The strict convexity of $\kappa_2$ thus implies that along the ray $x_1/x_2 = \gamma_2/\gamma_1$ it holds that $-\kappa_2'(-\gamma_2) < x_2/T(x_1, x_2) < -\kappa_2'(-\gamma_1)$. It is a matter of algebra to verify that these inequalities are equivalent to $s_2 < \gamma_2/\gamma_1 < s_1$ (see also Lemma 5 below). The assertions (i), (ii) and (iii) follow then in view of the definitions of $\mathcal{D}_i, i = 0, 1, 2$, and $\widehat{\mathcal{D}}_2$. □

A key point in the analysis of the asymptotics of the degenerate risk processes is an equivalent description of the cones $\mathcal{D}_i$ in terms of comparisons with the time $T = T(x_1, x_2)$ defined in (26), which will enable us to translate the asymptotics in two-dimensional space into the "space-time" asymptotics of Arfwedson (1955) and Höglund (1990).

LEMMA 5. *Writing $T_i = x_i/[-\kappa_i'(-\gamma_i)]$ and $\widetilde{T}_i = x_i/[-\kappa_i'(-\gamma_{3-i})]$, $i = 1, 2$, and $a(x) = \frac{x_1}{x_2}$ [where $x = (x_1, x_2) \in \mathbb{R}_+^2$], the following hold true:*

$$\mathcal{D}_2 = \{x \in \mathbb{R}_+^2 : a(x) < s_2\} = \{x \in \mathbb{R}_+^2 : \widetilde{T}_1 < T\} = \{x \in \mathbb{R}_+^2 : T_2 < T\},$$
$$\mathcal{D}_0 = \{x \in \mathbb{R}_+^2 : T_1 < T < T_2\},$$
$$\mathcal{D}_1 = \{x \in \mathbb{R}_+^2 : s_1 < a(x)\} = \{x \in \mathbb{R}_+^2 : T < \widetilde{T}_2\} = \{x \in \mathbb{R}_+^2 : T < T_1\}.$$

PROOF. The first two equalities for $\mathcal{D}_1$ and $\mathcal{D}_2$ are just the definitions, given in (47). Next, we note that from the definition of $T$ it is easy to check that

(49)     $x_1/x_2 = a \iff x_2/T(x_1, x_2) = (p_1 - p_2)/(1 - a) = v_a.$

In particular, inserting $a = s_1$ [defined in (48)] and using $\kappa_1'(s) = p_1 - p_2 + \kappa_2'(s)$ it follows that $s_1 = x_1/x_2$ iff $T = x_1/[-\kappa_1'(-\gamma_1)]$ iff $T = x_2/[-\kappa_2'(-\gamma_1)]$. The two representations for $\mathcal{D}_1$ now directly follow. The identities for $\mathcal{D}_2$ are proved similarly. Finally, the equalities for $\mathcal{D}_0$ follow by intersecting those for the complements $\overline{\mathcal{D}}_2^c$ and $\overline{\mathcal{D}}_1^c$. □

5.2. *Asymptotics.* The leading term asymptotics of the two-dimensional ruin probabilities will be expressed in terms of the usual "one dimensional large deviations cast": the adjustment coefficients $\gamma_i > 0$ of $X_i$ satisfying $\kappa_i(-\gamma_i) = 0$, and $\gamma(a)$ given for $0 < a < 1$ by

$$\gamma(a) = \kappa_2^*(-v_a)/v_a \qquad \text{where } v_a := (p_1 - p_2)/(1 - a).$$

Below we consider asymptotics along the rays $(aK, K)$ in the plane with $a < \overline{a}$, where

$$\overline{a} = 1 + (p_1 - p_2)/\underline{v}$$

with, as before,

(50)     $\underline{\theta} = \inf\{\theta : \kappa_2(\theta) < \infty\}$ and $\underline{v} = \lim_{\theta \downarrow \underline{\theta}} \kappa_2'(\theta).$

Note that in terms of the "space-time velocities" the restriction $a < \overline{a}$ reads as $v = v_a < -\underline{v}$.

THEOREM 5. *Assume that $\underline{\theta} < -\gamma_1$. If $a < \overline{a}$, it holds as $K \to \infty$,*

$$\psi_{\mathrm{sim}}(aK, K) \sim \begin{cases} C_1 e^{-\gamma_1 aK}, & \text{if } (aK, K) \in \mathcal{D}_1, \\ (D_2^\#(v_a) + D_2'(v_a))K^{-1/2}e^{-\gamma(a)K}, & \text{if } (aK, K) \in \mathcal{D}_0, \\ C_2 e^{-\gamma_2 K}, & \text{if } (aK, K) \in \mathcal{D}_2, \end{cases}$$

*where, for $i = 1, 2$, $C_i = -\kappa_i'(0)/\kappa_i'(-\gamma_i)$ and*

$$(51) \quad D_i'(w) = \frac{\theta_w^{(i)} - \theta_w}{|\theta_w \theta_w^{(i)}|} \frac{\sqrt{w}}{\sqrt{2\pi \kappa_{3-i}''(\theta_w)}},$$

$$(52) \quad D_i^\#(w) = \left[\frac{1}{\theta_w} - \frac{1}{\theta_w^{(i)}} + \frac{\kappa_{3-i}'(0)}{\kappa_{3-i}(\theta_w^{(i)})} - \frac{\kappa_{3-i}'(0)}{\kappa_{3-i}(\theta_w)}\right] \frac{\sqrt{w}}{\sqrt{2\pi \kappa_{3-i}''(\theta_w)}},$$

*where $\theta_w < \theta_w^{(i)}$ satisfy $\kappa_2'(\theta_w) = -w$ and $\kappa_i(\theta_w) = \kappa_i(\theta_w^{(i)})$.*

Next, we turn to the asymptotics of the ruin probability $\psi_{\mathrm{and}}$, which are formulated in terms of the (sometimes different) partition of the positive quadrant into $\mathcal{D}_1$, $\widehat{\mathcal{D}}_2$ and $\widehat{\mathcal{D}}_0 = \mathbb{R}_+^2 \setminus [\overline{\mathcal{D}}_1 \cup \overline{\widehat{\mathcal{D}}}_2]$.

THEOREM 6. *If $a < \overline{a}$, then it holds that, as $K \to \infty$,*

$$\psi_{\mathrm{and}}(aK, K) \sim \begin{cases} C_1 e^{-\gamma_1 aK}, & \text{if } (aK, K) \in \mathcal{D}_1, \\ (D_1'(v_a) - D_1^\#(v_a))K^{-1/2}e^{-\gamma(a)K}, & \text{if } (aK, K) \in \widehat{\mathcal{D}}_0, \\ \widehat{C}_2 e^{-(a\gamma_3 + (1-a)\gamma_2)K}, & \text{if } (aK, K) \in \widehat{\mathcal{D}}_2, \end{cases}$$

*where*

$$\widehat{C}_2 = -C_2 \kappa_1'(-\gamma_2)/\kappa_1'(-\gamma_3)$$

*and $D_1'$ and $D_1^\#$ are respectively given by (51) and (52).*

NOTE. These results imply that if $(aK, K)$ is contained in either $\mathcal{D}_1$ or $\mathcal{D}_2$, then $\psi_{\mathrm{sim}}(aK, K)$ and $\psi_{\mathrm{and}}(aK, K)$ are of the same order.

5.3. *Two terms asymptotic expansions in terms of one-dimensional shifted measures.* The key for obtaining "two terms asymptotic expansions" for $\psi_{\mathrm{or}}$ and $\psi_{\mathrm{sim}}$ (and then also leading term asymptotics) is given by the following decompositions that are generalizations of Corollary 1 to the current setting.

COROLLARY 2. *Let $x_2 > x_1$. It holds that*

$$\psi_{\mathrm{or}}(x_1, x_2) = \psi_1(x_1, T) + C_2(x_1, T)e^{-\gamma_2 x_2}\overline{\psi}_1^{(-\gamma_2)}(x_1, T),$$

$$\psi_{\mathrm{sim}}(x_1, x_2) = \psi_2(x_2, T) + C_1(x_2, T)e^{-\gamma_1 x_1}\overline{\psi}_2^{(-\gamma_1)}(x_2, T),$$

*where, for $i = 1, 2$,*

$$C_i(x_{3-i}, T) = E^{(-\gamma_i)}[h_i(X_{3-i}(T))|\tau_{3-i}(x_{3-i}) > T]$$

*with $h_i(x) = e^{\gamma_i x}\psi_i(x), i = 1, 2$.*

PROOF. In view of (30), the Markov property and a change of measure it follows that

$$(53) \quad \psi_{\text{or}}(x_1, x_2) = \psi_1(x_1, T) + P\Big(\tau_1(x_1) > T, \inf_{T \leq s < \infty} X_2(s) < 0\Big)$$

$$(54) \qquad\qquad = \psi_1(x_1, T) + e^{-\gamma_2 x_2} E^{(-\gamma_2)}[h_2(X_1(T))\mathbf{1}_{\{\tau_1(x_1)>T\}}]$$

using that $X_1(T) = X_2(T)$. The proof of the decomposition of $\psi_{\text{sim}}$ is similar and omitted. □

For $\psi_{\text{and}}$ similar decompositions are derived in the following result. Recall that $\gamma_3$ is defined as the largest root of $\kappa_1(-\gamma_2) = \kappa_1(-\theta)$ and set $\widetilde{\gamma} := \gamma_3 - \gamma_2$.

COROLLARY 3. *For $x_2 > x_1$ it holds that*

$$(55) \qquad \psi_{\text{and}}(x_1, x_2) = w_1(x_1, T) + P(\tau_1(x_1) \leq T, \tau_2(x_2) < \infty)$$

*where $w_1$ is given in (35) and, with $\tau_1 = \tau_1(x_1)$ and $\tau_2 = \tau_2(x_2)$,*

$$e^{\gamma_2 x_2 + \widetilde{\gamma} x_1} P(\tau_1 \leq T, \tau_2 < \infty)$$

$$(56) \qquad = E^{(-\widetilde{\gamma}, -\gamma_2)}[e^{\widetilde{\gamma} X_1(\tau_1)} h_2(X_2(\tau_1))\mathbf{1}_{\{\tau_1 \leq T\}}]$$

$$\qquad = e^{\widetilde{\gamma} x_1}\overline{C}_2(x_2, T)\psi_2^{(-\gamma_2)}(x_2, T) + \overline{C}_1(x_1, T)\psi_1^{(-\gamma_3)}(x_1, T),$$

*for $\overline{C}_2(x_2, T) = E^{(-\gamma_2)}[e^{\gamma_2 X_2(T)}\overline{\psi}_2(X_2(T))|\tau_2 \leq T]$ and*

$$\overline{C}_1(x_1, T) = E^{(-\gamma_3)}[e^{\gamma_3 X_1(T)}\psi_2(X_1(T))|\tau_1 \leq T].$$

PROOF. Recall that (55) was derived in Proposition 1. It follows by definition of $\gamma_3$ that $\kappa_1(-\gamma_2 - \widetilde{\gamma}) = \kappa_1(-\gamma_2)$ or, equivalently, $\kappa_1^{(-\gamma_2)}(-\widetilde{\gamma}) = 0$. In view of this observation and the form of $\kappa(u, v)$, derived in Example 1, it follows that $\kappa(-\widetilde{\gamma}, -\gamma_2) = 0$. Changing measure with the martingale $\exp(-\widetilde{\gamma}(X_1(t) - x_1) - \gamma_2(X_2(t) - x_2))$ and applying the strong Markov property at $\tau_1$ yields the first equality in (56).

The second equality follows by noting that

$$P(\tau_1 \leq T, \tau_2 < \infty)$$
$$\qquad = P(\tau_1 \leq T, \tau_2 \leq T) + P(\tau_1 \leq T, T < \tau_2 < \infty)$$

$$= P(\tau_2 \leq T) + E[\mathbf{1}_{\{\tau_1 \leq T < \tau_2\}} P_{X_2(T)}(\tau_2 < \infty)]$$
$$= E[\mathbf{1}_{\{\tau_2 \leq T\}} \overline{\psi}_2(X_2(T))] + E[\mathbf{1}_{\{\tau_1 \leq T\}} \psi_2(X_2(T))]$$
$$= e^{-\gamma_2 x_2} E^{(-\gamma_2)}[e^{\gamma_2 X_2(T)} \overline{\psi}_2(X_2(T))|\tau_2 \leq T] \psi_2^{(-\gamma_2)}(x_2, T)$$
$$+ e^{-\gamma_2 x_2 - \widetilde{\gamma} x_1} \times E^{(-\gamma_3)}[e^{\gamma_3 X_1(T)} \psi_2(X_1(T))|\tau_1 \leq T] \psi_1^{(-\gamma_3)}(x_1, T),$$

where in the second line we used the Markov property and that $\{\tau_2 \leq T\} \subset \{\tau_1 \leq T\}$ and in the last line we changed the measure and that $X_1(T) = X_2(T)$. □

We write

$$f \approx g + h \qquad \text{as } x \to \infty \text{ iff } \lim_{x \to \infty}(f - g)/h(x) = \lim_{x \to \infty}(f - h)/g(x) = 1.$$

PROPOSITION 2.  *Assume that $\underline{\theta} < -\gamma_1$. For any $v < -\underline{v}$, it holds that, if $x_1, x_2 \to \infty$ such that $x_2/T(x_1, x_2) = v$,*

(57)  $$\psi_{\mathrm{or}}(x_1, x_2) \approx \psi_1(x_1, T) + \widetilde{C}_2(v) e^{-\gamma_2 x_2} \overline{\psi}_1^{(-\gamma_2)}(x_1, T),$$

(58)  $$\psi_{\mathrm{sim}}(x_1, x_2) \approx \psi_2(x_2, T) + \widetilde{C}_1(v) e^{-\gamma_1 x_1} \overline{\psi}_2^{(-\gamma_1)}(x_2, T),$$

*where, for $i = 1, 2$ and $v \neq -\kappa_2'(-\gamma_1), -\kappa_2'(-\gamma_2)$,*

$$\widetilde{C}_i(v) = \begin{cases} C_i, & \text{if } -\kappa_2'(-\gamma_i) < v, \\ c_{3-i}(v, \gamma_i)^{-1}[\psi_i^*(\theta_v) - \psi_i^*(\theta_v^{(3-i)})], & \text{if } 0 < v < -\kappa_2'(-\gamma_i), \end{cases}$$

*where $\theta_v < \theta_v^{(i)}$ satisfy $\kappa_2'(\theta_v) = -v$, $\kappa_i(\theta_v) = \kappa_i(\theta_v^{(i)})$, and*

$$c_i(v, c) = \frac{\theta_v^{(i)} - \theta_v}{(\theta_v^{(i)} + c)(\theta_v + c)}$$

*with $\psi_i^*$ being the Laplace transform of $\psi_i$.*

PROPOSITION 3.  *For any $v < -\underline{v}$, it holds that, if $x_1, x_2 \to \infty$ such that $x_2/T(x_1, x_2) = v$,*

(59)  $$\begin{aligned} &\psi_{\mathrm{and}}(x_1, x_2) \\ &\approx \psi_1(x_1) - \psi_1(x_1, T) \\ &\quad + \{\overline{C}_2(v) e^{-\gamma_2 x_2} \psi_2^{(-\gamma_2)}(x_2, T) + \overline{C}_1(v) e^{-\gamma_2 x_2 - \widetilde{\gamma} x_1} \psi_1^{(-\gamma_3)}(x_1, T)\}, \end{aligned}$$

*where* $\widetilde{\gamma} = \gamma_3 - \gamma_2$ *and for* $v \neq -\kappa_2'(-\gamma_3)$,

$$\overline{C}_2(v) = \begin{cases} 0, & \text{if } 0 < v < -\kappa_2'(-\gamma_3), \\ |c_2(v, \gamma_2)|^{-1} \cdot \overline{\psi}_2^*(\theta_v^{(2)}), & \text{if } v > -\kappa_2'(-\gamma_3), \end{cases}$$

$$\overline{C}_1(v) = \begin{cases} -C_2 \dfrac{\kappa_1'(-\gamma_2)}{\kappa_1'(-\gamma_3)}, & \text{if } 0 < v < -\kappa_2'(-\gamma_3), \\ |c_1(v, \gamma_3)|^{-1} \cdot [\psi_2^*(\theta_v^{(1)}) - \theta_v^{-1}] & \text{if } v > -\kappa_2'(-\gamma_3), \end{cases}$$

*with* $\overline{\psi}_i^*$ *being the Laplace transform of* $\overline{\psi}_i$.

PROOF OF PROPOSITION 2. In view of Corollary 2, the proof for $\psi_{\text{or}}$ is complete once we show that $C_2(x_1, T)$ converges to $\widetilde{C}_2(v)$ if $x_1, x_2 \to \infty$ such that $x_2/T = v$. We distinguish between two cases.

If $x_1 + \kappa_1'(-\gamma_2)T > 0$, then the strong law of large numbers implies that $X_1(T) \to \infty$ $P^{(-\gamma_2)}$-a.s. and that $P^{(-\gamma_2)}(\tau_1(x_1) \leq T)$ tends to zero (see Theorem 1). Since $h_2(y) \to C_2$ as $y \to \infty$ (by the Cramér–Lundberg approximation), we conclude that $C_2(x_1, T)$ converges to $C_2$ (by bounded convergence).

In the case that $x_1 + \kappa_1'(-\gamma_2)T < 0$, we note that by virtue of Theorem 2(i) the distribution of $X_1(T)$ conditioned on $\tau_1 > T$ (under $P^{(-\gamma_2)}$) converges to the measure $\overline{\Psi}_v$ in (23), with the shifts in (23) calculated using the cumulant exponent $\kappa_1^{(-\gamma_2)}$ and direction $v' = x_1/T(x_1, x_2)$. Note that the application of Theorem 2 is justified since we have that $\underline{\theta} < -\gamma_1 < -\gamma_2$. Indeed, observe at this point that $\kappa_1^{(-\gamma_2)\prime}(\theta_v^{(-\gamma_2)}) = \kappa_1'(\theta_v) = -v'$ and that, in view of Remark 1, $\theta_v^{(-\gamma_2)} = \theta_v + \gamma_2$ and $\theta_v^{(1)(-\gamma_2)} = \theta_v^{(1)} + \gamma_2$. Thus, $C_2(x_2, T)$ converges to

(60) $$\int_0^\infty h_2(y)\overline{\Psi}_v(dy) = \int_0^\infty \psi_2(y) \frac{(\theta_v + \gamma_2)(\theta_v^{(1)} + \gamma_2)}{\theta_v^{(1)} - \theta_v} [e^{-\theta_v y} - e^{-\theta_v^{(1)} y}] \, dy.$$

The proof of the asymptotics of $\psi_{\text{sim}}$ is similar and omitted. □

PROOF OF PROPOSITION 3. In view of Corollary 3 to finish the proof, we have to show convergence of the conditional expectations in the two different cases.

In the first case when $x_1 + \kappa_1'(-\gamma_3)T < 0$, it follows by the law of large numbers that $P^{(-\gamma_3)}(\tau_1 \leq T)$ tends to 1. Also, taking note of Lemma 1 and of the fact that in view of the definition of $T$, it holds that $X_2(\tau_1) = X_1(\tau_1) + (p_1 - p_2)[T - \tau_1])$, it follows that $X_2(\tau_1) \to \infty$ and $h_2(X_2(\tau_1)) \to C_2$, $P^{(-\gamma_3)}$-a.s. Therefore, the bounded convergence theorem implies that the expectation in the first line of (56) converges to $C_2\widetilde{C}$ (where $\widetilde{C}$ denotes the asymptotic constant for $\psi_1$ under $P^{(-\gamma_2)}$).

In the opposite case that $x_1 + \kappa_1'(-\gamma_3)T > 0$ invoking Theorem 2 as in Proposition 2 (which is in this case justified as $\underline{\theta} \leq -\gamma_1 < -\gamma_3 \leq -\gamma_2$) yields the form of $\overline{C}_1(v)$ and $\overline{C}_2(v)$ and in view of (31), the proof is complete. □

5.4. *Proofs of Theorems 5 and 6.* In the proof, we use the following result:

LEMMA 6. (i) *We have* $\gamma(a) > \max\{a\gamma_1, \gamma_2\}$ *for* $a \neq s_1, s_2$ *and*

$$\gamma(a) = a\gamma_1 + \kappa_2^{*(-\gamma_1)}(-v_a)/v_a = a\gamma_1 + \kappa_1^{*(-\gamma_1)}(-av_a)/v_a = \kappa_1^*(-av_a)/v_a.$$

(ii) *If* $\kappa_i'(0^+) > 0$, *it holds that* $\psi_i^*(\theta) = \theta^{-1} - \kappa_i'(0)/\kappa_i(\theta)$.

PROOF. (i) is a direct consequence of the definitions of $\gamma(a), \gamma_i$ and Remark 1. (ii) directly follows from Bertoin (1996), Theorem VII.10.  □

PROOF OF THEOREMS 5 AND 6. In view of (58), the proof of Theorem 5 consists in identifying the leading order term by applying the Arfwedson–Höglund's Theorem 1 to the different terms in (58). Invoking the characterization of the cones in Lemma 5 and (49), it follows, for example, that for $(x_1, x_2) \in \mathcal{D}_2$, $\psi_2(x_2, T) \sim C_2 e^{-\gamma_2 x_2}$ whereas for $(x_1, x_2) \in \overline{\mathcal{D}}_1^c \supset \mathcal{D}_2$ it holds that

$$(61) \qquad \overline{\psi}_2^{(-\gamma_1)}(x_2, T) \sim |D_2^{(-\gamma_1)}(v_a)|(v_a/x_2)^{1/2} e^{-x_2 \kappa_2^{*(-\gamma_1)}(-v_a)/v_a},$$

where $D_2^{(-\gamma_1)}$ is specified by $D$ in Theorem 1 with $\kappa = \kappa_2^{(-\gamma_1)}$. Thus, by Lemma 6(i), it follows that the leading term is $C_2 e^{-\gamma_2 x_2}$ if $(x_1, x_2) \in \mathcal{D}_2$. Similarly, it can be checked that the leading term is $C_1 e^{-\gamma_1 x_1}$ if $(x_1, x_2) \in \mathcal{D}_1$. Finally, in the case that $(x_1, x_2) \in \mathcal{D}_0$, we note that both terms in (58) are of the same order [cf. Lemma 6(ii)]. More precisely,

$$\psi_{\mathrm{sim}}(x_1, x_2) \sim [D_2(v_a) + \widetilde{C}_1(v_a) D_2^{(-\gamma_1)}(v_a)](v_a/x_2)^{1/2} e^{-\gamma(a) x_2},$$

where $D_2, D_2^{(-\gamma_1)}$ are specified by $D$ in Theorem 1 with $\kappa = \kappa_2$ and $\kappa = \kappa_2^{(-\gamma_1)}$, respectively, and $\widetilde{C}_1$ is given in Proposition 2. Using Lemma 6 and that $\kappa_1'' = \kappa_2''$, it is a matter of algebra to verify the form of the constants $D_2'$ and $D_2^{\#}$.

Drawing on Proposition 3, Theorem 6 can be proved following an analogous line of reasoning. We omit the details.  □

**6. Examples.** We now develop two explicit examples that illustrate the results shown in the previous sections.

6.1. *Cramér–Lundberg model with exponential jumps.* Let $X$ be a drift $p$ minus a compound Poisson process with rate $\lambda$ and exponential jump sizes with mean $\mu$ starting at $x$. Then the characteristic function of $X$ reads as $\kappa(\theta) = p\theta - \lambda\theta/(\mu + \theta)$ and, if $p > \frac{\lambda}{\mu}$, the ultimate ruin probability is equal to $\psi(x) = Ce^{-\gamma x}$, where the adjustment coefficient is $\gamma = \mu - \lambda/p$ and $C =$ $\lambda/(\mu p)$. More generally, it was shown by Asmussen (1984), Knessl and Peters (1994) (with $p=1$) and Pervozvansky (1998) that the finite time ruin probability $\psi(x,t)$ is given by

$$\overline{\psi}(x,t) = 1 - \psi(x,t) = [1 - Ce^{-\gamma x}]\mathbf{1}_{(\gamma>0)} + w(x,t), \tag{62}$$

where

$$w(x,t) = \frac{1}{\pi}\sqrt{\frac{\lambda}{\mu p}} \int_{s_-}^{s_+} e^{a(q)x - qt} \sin(b(q)x - \phi(q)) \frac{dq}{q} \tag{63}$$

with $s_\pm = (\sqrt{\lambda} \pm \sqrt{\mu p})^2$, $\phi(q) = \arccos(\frac{p\mu + \lambda - q}{2\sqrt{\lambda\mu p}})$ and

$$a(q) = \frac{\lambda - \mu p - q}{2p}, \qquad b(q) = \frac{\sqrt{4pq\mu - (\lambda - \mu p - q)^2}}{2p}. \tag{64}$$

Further, we note that, under $P^{(c)}$, $X$ is still a drift $p$ minus a compound Poisson process with exponential jumps with the changed rates $\lambda_c = \lambda\frac{\mu}{\mu+c}$ and $\mu_c = \mu + c$. In particular, $\lambda_{-\gamma} = \mu p$ and $\mu_{-\gamma} = \lambda/p$ are the parameters under $P^{(-\gamma)}$.

In view of the previous paragraph, we see that, under $P^{(-\gamma_1)}$, the drift of $X_2$ is always negative, $\kappa_2^{(-\gamma_1)\prime}(0) = \kappa_2'(-\gamma_1) < 0$. Also, under $P^{(-\gamma_2)}$, the adjustment parameter of $X_1$ is positive if and only if $\rho > \rho^* := p_2^2/p_1$ and is then equal to

$$\widetilde{\gamma} = \gamma_3 - \gamma_2 = \frac{\mu}{p_2}\left(\rho - \frac{p_2^2}{p_1}\right), \tag{65}$$

and the asymptotic constant $\widetilde{C}$ satisfies $\widehat{C}_2 = \widetilde{C}C_2 = \frac{p_2}{p_1}$. Inserting the expressions (62)–(64) (with the proper choices of parameters) into Corollary 1 leads then to explicit expressions for $\psi_{\text{and}}, \psi_{\text{sim}}$ and $\psi_{\text{or}}$.

It is a matter of calculus to verify that

$$s_1 = \frac{p_1^2/\rho - p_1}{p_1^2/\rho - p_2}, \qquad s_2 = \frac{(p_2^2/\rho - p_1)_+}{p_2^2/\rho - p_2}$$

and, if $\rho > \rho^*$,

$$s_3 = \frac{\rho p_1^2/p_2^2 - p_1}{\rho p_1^2/p_2^2 - p_2}.$$

Also, by invoking Corollary 1 or by a direct calculation, we see that

$$C_i(v) \equiv \frac{\lambda}{\mu p_i} = C_i, \qquad i = 1, 2.$$

Inserting these quantities into Propositions 2 and 3 yields explicit asymptotics expansions for $\psi_{\text{and}}, \psi_{\text{sim}}$ and $\psi_{\text{or}}$.

6.2. *Brownian motion with drift.* If $X(t) = mt + B(t)$ where $B(t)$ is standard Brownian motion, then its characteristic exponent reads as $\kappa(\theta) = \frac{1}{2}\theta^2 + m\theta$. If $m > 0$, $\psi(x) = e^{-\gamma x}$, where $\gamma = 2m$ is the adjustment coefficient. Further, under $P^{(c)}$, $X$ is still Brownian motion, but the drift changes to $m + c$. The drift of the measure associated to $c = -\gamma$ is $-m$, that is, the Brownian motion switches its drift. In view of the Corollary 1 and the well-known first-passage distribution of Brownian motion with drift,

$$\overline{\psi}(x,t) = \Phi\left(\frac{x+mt}{\sqrt{t}}\right) - e^{-2mx}\Phi\left(\frac{-x+mt}{\sqrt{t}}\right), \tag{66}$$

we find that if $x_2 > x_1$ then

$$\begin{aligned}
\psi_{\text{or}}(x_1, x_2) &= P(\tau_1(x_1) \leq T) + e^{-2p_2 x_2} P^{(-2p_2)}(\tau_1(x_1) > T) \\
&= 1 - \Phi(a(x_1, p_1)) + e^{-2p_1 x_1}\Phi(a(-x_1, p_1)) \\
&\quad + e^{-2p_2 x_2} \times [\Phi(a(x_1, p_1 - 2p_2)) \\
&\qquad - e^{-2x_1(p_1 - 2p_2)}\Phi(a(-x_1, p_1 - 2p_2))],
\end{aligned} \tag{67}$$

$$\begin{aligned}
\psi_{\text{sim}}(x_1, x_2) &= P(\tau_2(x_2) \leq T) + e^{-2p_1 x_1} P^{(-2p_1)}(\tau_2(x_1) > T) \\
&= 1 - \Phi(a(x_2, p_2)) + e^{-2p_2 x_2}\Phi(a(-x_2, p_2)) \\
&\quad + e^{-2p_1 x_1} \times [\Phi(a(x_2, p_2 - 2p_1)) \\
&\qquad - e^{-2x_2(p_2 - 2p_1)}\Phi(a(-x_2, p_2 - 2p_1))],
\end{aligned} \tag{68}$$

$$\begin{aligned}
\psi_{\text{and}}(x_1, x_2) &= \Phi(a(x_1, p_1)) - 1 + e^{-2p_1 x_1}(1 - \Phi(a(-x_1, p_1))) \\
&\quad + e^{-2p_2 x_2} \times [e^{-2x_1(p_1 - 2p_2)}\Phi(a(-x_1, p_1 - 2p_2)) \\
&\qquad + 1 - \Phi(a(x_1, p_1 - 2p_2))],
\end{aligned} \tag{69}$$

where $a(x, p) = [x + pT]/\sqrt{T}$ and $\Phi$ denotes the cumulative standard normal distribution function. In view of the facts that $\Phi(-x) = 1 - \Phi(x)$ and $1 - \Phi(x) \sim (2\pi)^{-1/2} x^{-1} \exp(-x^2/2)$ as $x \to \infty$, it follows from (68) and (69) that if $x_1, x_2$ tend to infinity with $x_1/x_2 = a$ then

$$\psi_{\text{and}}(x_1, x_2) \approx \begin{cases} e^{-2p_2 x_2 - 2(p_1 - 2p_2)^+ x_1} + o(av_a), & \text{if } 0 < a < s_3, \\ o(av_a), & \text{if } s_3 < a < s_1, \\ e^{-2p_1 x_1} + o(av_a), & \text{if } s_1 < a < 1, \end{cases} \tag{70}$$

$$\psi_{\text{sim}}(x_1, x_2) \approx \begin{cases} e^{-2p_2 x_2} + \delta(v_a), & \text{if } 0 < a < s_2, \\ \delta(v_a), & \text{if } s_2 < a < s_1, \\ e^{-2p_1 x_1} + \delta(v_a), & \text{if } s_1 < a < 1, \end{cases} \tag{71}$$

where

$$o(v) = \left[\frac{2v}{p_1^2 - v^2} + \frac{2v}{v^2 - (p_1 - 2p_2)^2}\right] \frac{\sqrt{v}}{\sqrt{2\pi x_1}} e^{-x_1(v + p_1)^2/[2v]},$$

$$\delta(v) = \left[\frac{2v}{v^2 - p_2^2} + \frac{2v}{(p_2 - 2p_1)^2 - v^2}\right] \frac{\sqrt{v}}{\sqrt{2\pi x_2}} e^{-x_2(v+p_2)^2/[2v]}$$

and

$$s_1 = \frac{p_1}{2p_1 - p_2}, \qquad s_2 = \frac{(2p_2 - p_1)_+}{p_2},$$

and, if $p_1 > 2p_2$,

$$s_3 = \frac{p_1 - 2p_2}{2p_1 - 3p_2}.$$

The asymptotics of $\psi_{\text{sim}}$ agree with the asymptotics of the steady state distribution of a tandem queue calculated in Lieshout and Mandjes (2007).

By straightforward calculations, it can be verified that if $S$ is a Brownian motion, then $C_i = 1$ and $\theta_v = -v - p_2$, $\theta'_v = v - p_2$, $\theta_v^\star = v + p_2 - 2p_1$ and

$$\widetilde{C}_i(v) = 1, \qquad \kappa_i^*(-v) = \frac{(v + p_i)^2}{2}, \qquad \kappa_1^*(-av_a) = \kappa_2^*(-v_a)$$

for $i = 1, 2$. Inserting these quantities in (58)–(59) and comparing with (70) and (71) it follows that Propositions 2 and 3 and Theorems 5 and 6 remain valid if $S$ is a Brownian motion.

## APPENDIX

In this section, we briefly review the general large deviations theory for first passage times developed by Collamore (1996) and explicitly relate it to the results for $\psi_{\text{sim}}$ in Theorem 5 by calculating the relevant quantities of the general theory.

We denote by $(X_1, X_2)$ the two-dimensional Lévy process given by $X_i(t) = p_i t - S(t)$ where $S$ is a Lévy process (starting at zero) and write $A = A_a = (-\infty, -a) \times (-\infty, -1)$ and $T_B$ for the first hitting time of a set $B$ by $(X_1, X_2)$. The probability $\psi_{\text{sim}}(aK, K)$ is then equal to $P(T_{KA} < \infty)$. If $\kappa(\theta_1, \theta_2) = \log E[e^{\theta_1 X_1 + \theta_2 X_2}]$ is finite in a neighborhood of the Cramér set $\mathcal{C} = \{\theta : \kappa(\theta) \leq 0\}$ and $0 \in \mathcal{C}^o$, then it holds that [Theorem 2.1 of Collamore (1996)]

$$\lim_{K \to \infty} \frac{1}{K} \log \psi_{\text{sim}}(aK, K) = -\inf_{x \in A_a} \widetilde{I}(x),$$

where $\widetilde{I}$ is the support function of the Cramér set $\mathcal{C}$,

(72) $$\widetilde{I}(x) = \sup_{\theta \in \mathcal{C}} \langle \theta, x \rangle.$$

In the next result $\widetilde{I}(x)$ is identified:

PROPOSITION 4. *For $x_1, x_2 < 0$, it holds that*
$$\widetilde{I}(x_1, x_2) = |x_2|\gamma(x_1/x_2),$$
*where*
$$\gamma(a) = \begin{cases} \kappa_2^*(-|v_a|)/|v_a|, & \text{if } a \neq 1, \\ -\underline{\theta} := -\inf_{\theta \in \mathcal{C}} \langle \theta, \mathbf{1} \rangle, & \text{if } a = 1 \end{cases}$$
*with $\mathbf{1} = (1,1)$, $\kappa_2^*(s) := \sup_{\theta \in \mathbb{R}} \{\theta s - \kappa_2(\theta)\}$ and $v_a = (p_1 - p_2)/(1-a)$.*

Using this result, we can calculate the first-passage rate function:

COROLLARY 4. *It holds for $a > 0$ that*
$$\inf_{x \in A_a} \widetilde{I}(x) = \begin{cases} \gamma_2, & \text{if } 0 < a \leq s_2, \\ \gamma(a), & \text{if } s_2 < a < s_1, \\ a\gamma_1, & \text{if } a \geq s_1. \end{cases}$$

We observe here that Corollary 4 agrees with the exponents of the asymptotics found in Theorem 5 (as it should).

PROOF OF PROPOSITION 4. The linear functional $\theta \mapsto \langle \theta, x \rangle$ attains it maximum over the closed set $\mathcal{C}$ at a point $\theta^*$ of the boundary $\partial \mathcal{C}$. In view of the form of $\kappa$ (Example 1), the definition of $v_a$ and $\kappa(\theta^*) = 0$, it is a matter of algebra to check that for $a \neq 1$
$$\widetilde{I}(x_1, x_2) = \sup_{\theta \in \partial \mathcal{C}} \left\{ \theta_1 x_1 + \theta_2 x_2 - (\kappa_2(\theta_1 + \theta_2) + (p_1 - p_2)\theta_1) \cdot \left|\frac{x_2}{v_a}\right| \right\}$$
$$= \left|\frac{x_2}{v_a}\right| \cdot \sup_{\theta \in \partial \mathcal{C}} \{(\theta_1 + \theta_2)(-|v_a|) - \kappa_2(\theta_1 + \theta_2)\}$$
$$= \left|\frac{x_2}{v_a}\right| \sup_{\eta \geq \underline{\theta}} \{\eta(-|v_a|) - \kappa_2(\eta)\} = \left|\frac{x_2}{v_a}\right| \kappa_2^*(-|v_a|),$$
where $\underline{\theta} = \inf\{\theta : \kappa_2(\theta) < \infty\}$ and in the last line we used that $\kappa_2(\theta_1 + \theta_2) = -\theta_1(p_1 - p_2)$ for $(\theta_1, \theta_2) \in \partial \mathcal{C}$. The rest of the statements follows by straightforward calculations. □

PROOF OF COROLLARY 4. From the form of $\widetilde{I}$, we deduce that $\widetilde{I}$ attains it minimum over $A_a$ at the boundary $\partial A_a$. If we set $x_2 = -1$, the minimization reduces to $\inf_{b \geq a} \gamma(b)$. Taking note of the fact that $\kappa_2^*(-v) \geq \gamma_2 v$ with equality if and only if $v = -s_2$ we see that

(73) $$\inf_{b \geq a} \gamma(b) = \begin{cases} \gamma_2, & \text{if } 0 < a \leq s_2, \\ \gamma(a), & \text{if } a > s_2. \end{cases}$$

Similarly, setting $x_1 = -a$, leads to the minimization $\inf_{c \leq 1} c\gamma(a/c)$ or equivalently, $\inf_{d \leq a} \gamma(d)/d$. Observing that

$$\frac{\kappa_2^*(-|v_a|)}{a|v_a|} = \sup_{\theta \in \mathbb{R}} \left\{ -\theta - \frac{\kappa_1(\theta)}{a|v_a|} \right\},$$

we see that $\kappa_2^*(-|v_a|) \geq a\gamma_1|v_a|$ with equality if and only if $a = s_1$ and conclude that

(74) $$\inf_{d \leq a} \gamma(d)/d = \begin{cases} a\gamma_1, & \text{if } a \geq s_1, \\ \gamma(a), & \text{if } 0 < a < s_1. \end{cases}$$

Combining equations (73) and (74) completes the proofs. □

Let us also note that classical arguments for the large deviations (LD) theory for stationary increments processes allow us to explain heuristically the structure of "or" and "sim" ruins in the quadrant. Let us consider for example the "sim" ruin. From the LD theory, we expect that paths exiting the positive quadrant, seen from far away, will be concentrated near one of three possible directions: the direct path to the origin, and the "dominant" (most probable) paths reaching the $x_1 = 0/x_2 = 0$ axes, respectively.

It turns out that all the dominant exit paths for reaching the $\{x_1 < 0, x_2 = 0\}$ semi-axis are parallel to each other and, therefore, give rise to a "boundary" cone $\mathcal{D}_1$. Similarly, the dominant exit paths reaching the $x_1 = 0, x_2 < 0$ semi-axis gives rise to a cone $\mathcal{D}_2$. These two cones will be disjoint in our case, as indicated in Figure 2. In these two boundary cones, the probability of simultaneous ruin is equivalent to the probability of ruin of the $X_1/X_2$ process respectively, that is, it holds asymptotically for large $x_1, x_2$ on a ray that lies within these cones $\mathcal{D}_i$ that $\psi_{\text{sim}}(x_1, x_2) \approx \psi_i(x_i)$.

A similar result holds for the "or" ruin. We may note geometrically that the boundary cones for hitting the semi-axes $\{x_1 = 0, x_2 > 0\}$ and $x_2 = 0, x_1 > 0$ are precisely the complements $\mathcal{C}_2 = \overline{\mathcal{D}}_1^c, \mathcal{C}_1 = \overline{\mathcal{D}}_2^c$.

This intuitive picture is confirmed and sharpened by Theorem 5 and Proposition 2.

**Acknowledgments.** We are much indebted to an anonymous referee for the careful reading of the manuscript and to Tomasz Rolski for his continuing advice on this paper.

F. Avram
Departement de Mathematiques
Université de Pau
64000 Pau
France
E-mail: Florin.Avram@univ-Pau.fr

Z. Palmowski
Department of Mathematics
University of Wrocław
pl. Grunwaldzki 2/4
50-384 Wrocław
Poland
E-mail: zpalma@math.uni.wroc.pl

M. R. Pistorius
Department of Mathematics
King's College London
The Strand
London WC2R 2LS
United Kingdom
E-mail: Martijn.Pistorius@kcl.ac.uk